\documentclass[11pt,draftcls,onecolumn]{IEEEtran}

\usepackage[sort]{cite}
\usepackage{psfrag}
\usepackage{epsfig}
\usepackage{amsmath,bbm,times,stmaryrd}

\usepackage{amsthm}

 \usepackage[mathscr]{eucal}

\usepackage{tabularx}
\usepackage{multirow}
\usepackage{enumerate}

\usepackage[normal]{subfigure}

\usepackage{colortbl}
\usepackage{dcolumn}
\newcolumntype{.}{D{.}{.}{1.3}}

\usepackage{txfonts}
\usepackage[subnum]{cases}

\usepackage{equationarray}

\newcommand\dashrule{\leavevmode\xleaders\hbox{-}\hfill\kern0pt}

\def\diag{\operatorname{diag}}

\newcommand{\ba}{{\bvec{a}}}

\newcommand{\bd}{\bvec{d}}

\newcommand{\bg}{\bvec{g}}

\newcommand{\bm}{\bvec{m}}

\newcommand{\bt}{\bvec{t}}
\newcommand{\bu}{\bvec{u}}
\newcommand{\bv}{\bvec{v}}
\newcommand{\bw}{\bvec{w}}

\newcommand{\bz}{\bvec{z}}

\newcommand{\bA}{{\bf A}}
\newcommand{\bB}{{\bf B}}
\newcommand{\bC}{{\bf C}}

\newcommand{\bF}{{\bf F}}
\newcommand{\bG}{{\bf G}}
\newcommand{\bH}{{\bf H}}
\newcommand{\bI}{{\bf I}}

\newcommand{\bK}{{\bf K}}

\newcommand{\bQ}{{\bf Q}}
\newcommand{\bR}{{\bf R}}

\newcommand{\bU}{{\bf U}}
\newcommand{\bV}{{\bf V}}
\newcommand{\bW}{{\bf W}}
\newcommand{\bX}{{\bf X}}
\newcommand{\bY}{{\bf Y}}
\newcommand{\bZ}{{\bf Z}}

\newcommand{\calO}{{\mathcal{O}}}

\newcommand{\calR}{{\mathcal{R}}}

\newcommand{\calX}{{\mathcal{X}}}

\newcommand{\bbeta}{\mbox{\boldmath $\beta$}}

\newcommand{\blambda}{\mbox{\boldmath $\lambda$}}

\newcommand{\bxi}{\mbox{\boldmath $\xi$}}

\newcommand{\bsigma}{\mbox{\boldmath $\sigma$}}

\newcommand{\bGamma}{\mbox{\boldmath $\Gamma$}}

\newcommand{\bSigma}{\mbox{\boldmath $\Sigma$}}

\newcommand{\bPhi}{\mbox{\boldmath $\Phi$}}
\newcommand{\bPsi}{\mbox{\boldmath $\Psi$}}

\newcommand{\Real}{\mathbb R}

\newcommand{\0}{\mbox{\boldmath $0$}}

\newcommand{\be}{\begin{eqnarray}}
\newcommand{\ee}{\end{eqnarray}}

\newcommand{\matrixb}{\left[ \begin{array}}
\newcommand{\matrixe}{\end{array} \right]}

\newcommand{\tr}{\mathop{\rm tr}\nolimits}

\def\*{\circledast}

\newtheorem{lemma}{Lemma}
\newtheorem{theorem}{{Theorem}}

\newcommand{\bvec}[1]{\boldsymbol{#1}}

\def\vectorize{\operatorname{vec}}
\newcommand{\vtr}[1]{\vectorize\hspace{-.3ex}\left(#1\right)}

\newcommand{\tensor}[1]{\boldsymbol{\mathscr{\MakeUppercase{#1}}}} 
\newcommand{\tA}{\tensor{A}}

\newcommand{\tE}{\tensor{E}}

\newcommand{\tG}{\tensor{G}}
\newcommand{\tH}{\tensor{H}}

\newcommand{\tS}{\tensor{S}}

\newcommand{\tX}{\tensor{X}}
\newcommand{\tY}{\tensor{Y}}
\newcommand{\tZ}{\tensor{Z}}

\usepackage[vlined,ruled,commentsnumbered]{algorithm2e}

\usepackage{aurical}

\renewcommand{\bigodot}{\mathop{\mbox{\fontsize{18}{19}\selectfont$\odot$}}}
\renewcommand{\bigotimes}{\mathop{\mbox{\fontsize{18}{19}\selectfont$\otimes$}}}

\usepackage{arydshln}

\usepackage{url}
\usepackage{soul,color}
\setstcolor{red}

\newcommand{\eqdef}{\stackrel{\triangle}{=}}

\setstcolor{red}

\renewcommand{\CommentSty}[1]{\fontsize{8.7}{9.8}\selectfont\textnormal{\texttt{#1}}\unskip}
\SetKwComment{mtcc}{\% }{}
\SetKwSwitch{Switch}{Case}{Other}{switch}{}{case}{otherwise}{end switch}%

\newcounter{example}

\newenvironment{example}
{\refstepcounter{example}\vspace{10pt}\par\noindent 
\textbf{Example \theexample\ }
\begin{itshape}}
{\end{itshape}}%

\usepackage[usenames,dvipsnames,svgnames,table]{xcolor}
\usepackage{array}

\usepackage{tikz}

\newcommand{\bcol}{Emerald}

\makeatletter
\tikzset{%
     remember picture with id/.style={%
       remember picture,
       overlay,
       draw=\bcol,
       save picture id=#1,
     },
     save picture id/.code={%
       \edef\pgf@temp{#1}%
       \immediate\write\pgfutil@auxout{%
         \noexpand\savepointas{\pgf@temp}{\pgfpictureid}}%
     },
     if picture id/.code args={#1#2#3}{%
       \@ifundefined{save@pt@#1}{%
         \pgfkeysalso{#3}%
       }{
         \pgfkeysalso{#2}%
       }
     }
   }

   \def\savepointas#1#2{%
  \expandafter\gdef\csname save@pt@#1\endcsname{#2}%
}

\def\tmk@labeldef#1,#2\@nil{%
  \def\tmk@label{#1}%
  \def\tmk@def{#2}%
}

\tikzdeclarecoordinatesystem{pic}{%
  \pgfutil@in@,{#1}%
  \ifpgfutil@in@%
    \tmk@labeldef#1\@nil
  \else
    \tmk@labeldef#1,(0pt,0pt)\@nil
  \fi
  \@ifundefined{save@pt@\tmk@label}{%
    \tikz@scan@one@point\pgfutil@firstofone\tmk@def
  }{%
  \pgfsys@getposition{\csname save@pt@\tmk@label\endcsname}\save@orig@pic%
  \pgfsys@getposition{\pgfpictureid}\save@this@pic%
  \pgf@process{\pgfpointorigin\save@this@pic}%
  \pgf@xa=\pgf@x
  \pgf@ya=\pgf@y
  \pgf@process{\pgfpointorigin\save@orig@pic}%
  \advance\pgf@x by -\pgf@xa
  \advance\pgf@y by -\pgf@ya
  }%
}
\makeatother

\graphicspath{{Fig_rev1/}}

\title{Tensor Deflation for CANDECOMP/PARAFAC. Part 3: Rank Splitting} 

\author{Anh-Huy~Phan$^{*}$, Petr~Tichavsk{\'y} and Andrzej~Cichocki
\thanks{A. H. Phan and A. Cichocki are with the Lab for Advanced Brain Signal Processing, Brain Science Institute, and BSI-TOYOTA Collaboration Center, RIKEN, Wakoshi, Japan, e-mail: (phan,cia)@brain.riken.jp.}
\thanks{A. Cichocki is also with System Research Institute, Warsaw, Poland.}
\thanks{P.  Tichavsk{\'y} is with Institute of Information Theory and Automation, Prague, Czech Republic, email: tichavsk@utia.cas.cz.}
\thanks{The work of P. Tichavsk{\'y} was supported by the Czech Science Foundation through project No. 14--13713S.}
}

\def\comment#1{}

\begin{document}

\maketitle

\begin{abstract}
\small
CANDECOMP/PARAFAC (CPD) approximates multiway data by sum of rank-1 tensors. 
Our recent study has presented a method to rank-1 tensor deflation, i.e. sequential extraction of the rank-1 components.
In this paper, we extend the method to block deflation problem. 
When at least two factor matrices have full column rank,  one can
extract two rank-1 tensors simultaneously, and rank of the data tensor is reduced by 2. For decomposition of order-3 tensors of size $R \times R \times R$ and rank-$R$, the block deflation has a complexity of $\calO(R^3)$ per iteration which is lower than the cost $\calO(R^4)$ of the ALS algorithm for the overall CPD.
\end{abstract}


\section{Introduction}

An important property in matrix factorisations like eigenvalue decomposition or singular value decomposition, is that rank-1 matrix components can be sequentially estimated via deflation method, such as the power iteration method.
The matrix deflation procedure is possible because subtracting the best rank-1 term from a matrix reduces the matrix rank. 
Unfortunately, this sequential extraction procedure in general is not applicable to decompose a rank-$R$ tensor \cite{Stegeman:2010:SBR}.


%
%

In our recent study\cite{Phan_ALS_deflation,Phan_tensordeflation_alg}, we have introduced a tensor decomposition which is able to extract a rank-1 tensor from a high rank tensor. The method is based on the rank-1 plus multilinear-$(R-1,R-1,R-1)$ block tensor decomposition, but with a smaller number of parameters, only two vectors per modes. This paper extends the rank-1 tensor extraction to block tensor deflation or rank splitting which splits a high rank-$R$ tensor into two tensors with smaller ranks. In particular, we develop an alternating subspace update (ASU) algorithm to extract a multilinear rank-(2,2,2) tensor from a rank-$R$ tensor. Since decomposition of a $2 \times 2 \times 2$ tensor can be found in closed-form, we can straightforwardly obtain the desired rank-1 components. The proposed algorithm estimates only 4 vectors and two scalars per dimension with a computational complexity of $\calO(R^3)$.
Moreover, it also requires a lower space cost than algorithms for the ordinary CANDECOMP/PARAFAC (CPD).

The paper is organised as follows.
A tensor decomposition for block tensor deflation or rank splitting is presented  in Section~\ref{sec::preliminaries}.
The proposed algorithm is presented in Section \ref{sec::asu_blk2}.
Simulations in Section~\ref{sec::simulations} will verify validity and performance of the proposed algorithm.
Section~\ref{sec::conclusion} concludes the paper.

\section{Preliminaries}\label{sec::preliminaries}

Throughout the paper, we shall denote tensors by bold calligraphic letters, e.g., $\tA \in \Real^{I_1 \times I_2 \times \cdots \times I_N}$,
matrices by bold capital letters, e.g., $\bA$ =$[\ba_1,\ba_2, \ldots, \ba_R] \in \Real^{I \times R}$, and vectors by bold italic letters, e.g., $\ba_j$.
%
The Kronecker product is denoted by  $\otimes$. Inner product of two tensors is denoted by $\langle \tX, \tY\rangle = \vtr{\tX}^T \vtr{\tY}$. Contraction between two tensors along modes-$\bm$, where $\bm = [m_1,\ldots,m_K]$, is denoted by $\langle \tX, \tY\rangle_{\bm}$, whereas $\langle \tX, \tY\rangle_{-n}$ represents contraction along all modes but mode-$n$.
Generally, we adopt notation used in \cite{NMF-book}.
%

The mode-$n$ matricization of tensor $\tY$ is denoted by $\bY_{(n)}$.
The mode-$n$ multiplication of a tensor
${\tY} \in \Real^{I_{1} \times I_{2} \times \cdots \times
I_{N}}$ by a matrix $\bU \in \Real^{I_n \times R}$ is denoted by
${\tZ} = {\tY} \;
 \times_n \; \bU \in \Real^{I_1 \times \cdots \times I_{n-1} \times R
\times  I_{n+1} \times \cdots \times I_N}$.
Products of a tensor $\tY$ with a set of $N$ matrices
$\{\bU^{(n)}\} = \left\{  \bU^{(1)},  \bU^{(2)}, \ldots,\right.$ $\left.\bU^{(N)}\right\}$
are denoted by
$
    {\tY}\, {\times} \, \{\bU^{(n)}\} \eqdef {\tY}  \,  {\times}_1  \, \bU^{(1)}
    \,  {\times}_2 \, \bU^{(2)} \cdots  {\times}_N  \, \bU^{(N)}$.

A tensor $\tX \in \Real^{I_1 \times I_2 \times \cdots \times I_N}$ is said in Kruskal form if  
\be
{\tX}   =   \sum\limits_{r = 1}^R  \lambda_r \, {\ba^{(1)}_{r}  \circ \ba^{(2)}_{r} \circ  \cdots  \circ \ba^{(N)}_{r}}   \, ,
\label{equ_CP_nolambda}
\ee
where ``$\circ$"  denotes the outer product,  $\bA^{(n)}=[\ba^{(n)}_1, \ba^{(n)}_2,\ldots,\ba^{(n)}_R]$ $ \in \Real ^{I_n \times R}$ are factor matrices,
$\ba^{(n) T}_r \ba^{(n)}_r = 1$, for $r = 1, \ldots, R$ and $n = 1, \ldots, N$, and $\lambda_1 \ge \lambda_2 \ge \cdots \ge \lambda_R >0$.


A tensor $\tX \in \Real^{I_1 \times I_2 \times \cdots \times I_N}$ has multilinear rank-$(R_1,R_2,\ldots,R_N)$ if $\emph{rank}(\bX_{(n)}) = R_n \le I_n$ for $n = 1, \ldots, N$, and can be expressed in the Tucker form as
\be
{\tX}  &=&  \sum\limits_{r_1 = 1}^{R_1}\sum\limits_{r_2 = 1}^{R_2}\cdots\sum\limits_{r_n = 1}^{R_N}  g_{r_1r_2\ldots r_N} \, {\ba^{(1)}_{r_1}  \circ \ba^{(2)}_{r_2} \circ  \cdots  \circ \ba^{(N)}_{r_N}}  ,  \label{equ_tucker}
\ee 
where $\tG = [g_{r_1r_2\ldots r_N}]$, and $\bA^{(n)}$ are of full column rank.
For compact expression, $\llbracket \blambda; \{\bA^{(n)}\} \rrbracket$ denotes a Kruskal tensor, where $\llbracket \bG; \{\bA^{(n)}\}\rrbracket$ represents a Tucker tensor.

\begin{figure}[t]
\centering
\psfrag{IxJxK}[c][c]{\scalebox{1}{\color[rgb]{0,0,0}\setlength{\tabcolsep}{0pt}\begin{tabular}{c}\footnotesize $I \times J \times K $\end{tabular}}}%

\psfrag{KxKxK}[t][t]{\scalebox{1}{\color[rgb]{0,0,0}\setlength{\tabcolsep}{0pt}\begin{tabular}{c}\\[-1.5em] \end{tabular}}}%
\psfrag{JxJxJ}[tl][tl]{\scalebox{1}{\color[rgb]{0,0,0}\setlength{\tabcolsep}{0pt}\begin{tabular}{l}\\[-1.5em] \end{tabular}}}


\psfrag{Y}[c][c]{\scalebox{1}{\color[rgb]{0,0,0}\setlength{\tabcolsep}{0pt}\begin{tabular}{c} \Large $\tY$\end{tabular}}}%
\psfrag{A}[c][c]{\scalebox{1}{\color[rgb]{0,0,0}\setlength{\tabcolsep}{0pt}\begin{tabular}{c} \Large $\bA$\end{tabular}}}%
\psfrag{B}[c][c]{\scalebox{1}{\color[rgb]{0,0,0}\setlength{\tabcolsep}{0pt}\begin{tabular}{c}\Large $\bB^T$\end{tabular}}}%
\psfrag{C}[c][c]{\scalebox{1}{\color[rgb]{0,0,0}\setlength{\tabcolsep}{0pt}\begin{tabular}{c}\Large $\bC$\end{tabular}}}%
\psfrag{G}[c][c]{\scalebox{1}{\color[rgb]{0,0,0}\setlength{\tabcolsep}{0pt}\begin{tabular}{c}  \Large $\tG$\end{tabular}}}%
\psfrag{H}[c][c]{\scalebox{1}{\color[rgb]{0,0,0}\setlength{\tabcolsep}{0pt}\begin{tabular}{c}  \Large $\tH$\end{tabular}}}%


\includegraphics[width=.8\linewidth, trim = 0.0cm .6cm 0cm 0cm,clip=false]{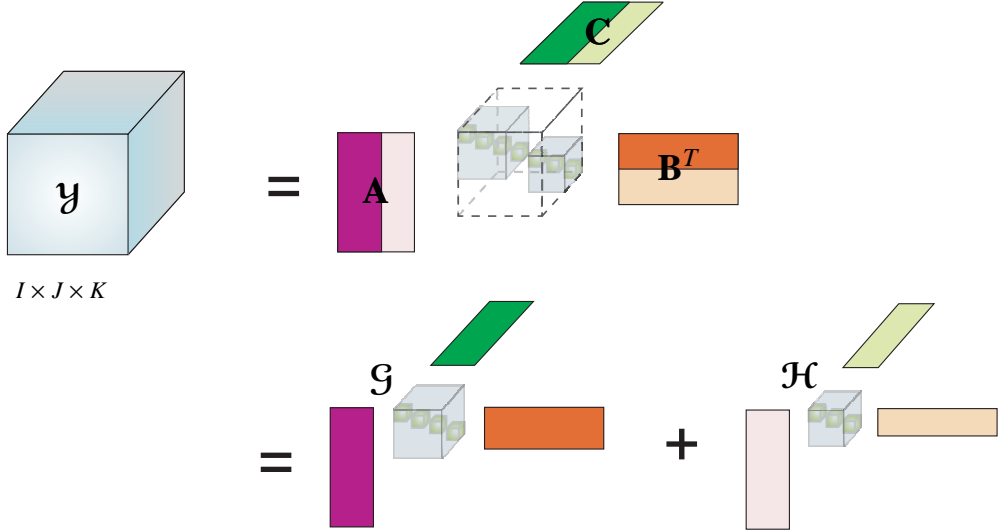}
\caption{Rank splitting for CP decomposition of a rank-$R$ tensor into two multilinear rank-$(K,\ldots,K)$ and rank-$(R-K,\ldots,R-K)$ tensors $\tG$ and $\tH$.
}
\label{cpd_rank_splitting}
\end{figure}

The main focus of this paper is
a block deflation which splits a rank-$R$ CPD into two sub rank-$K$  and rank-$(R-K)$ CPDs. 
This tensor decomposition is a particular case of the block tensor decomposition \cite{Lath-Nion-BCM3} but with only two blocks of multilinear rank-$(K,K,K)$ and rank-$(R-K,R-K,R-K)$ as illustrated in Fig.~\ref{cpd_rank_splitting}. 
{
{That is 
\be
\tY \approx \llbracket \tG ; \bU^{(1)}, \bU^{(2)},\ldots, \bU^{(N)}\rrbracket  + \llbracket \tH ; \bV^{(1)}, \bV^{(2)}, \ldots, \bV^{(N)}\rrbracket + \tE \label{equ_btd}
\ee
where $\bU^{(n)}$ and $\bV^{(n)}$ are matrices of size $I_n \times K$ and $I_n \times (R-K)$, respectively.}}
Following this tensor decomposition, decomposition of a rank-$R$ tensor can  proceed simultaneously through decompositions of sub-tensors with smaller ranks.
When $K = 1$, we have the rank-1 tensor deflation discussed in Part-1 \cite{Phan_tensordeflation_alg} and Part-2 \cite{Phan_tensordeflation_init}.

For this kind of tensor decomposition and block tensor deflation, we can use the ALS algorithm \cite{Lath-Nion-BCM3} or the non-linear least squares  (NLS) algorithm \cite{Sorber_fusion2013} developed for the multilinear rank-$(L_r,M_r,N_r)$ block tensor decomposition with two blocks. However, these existing algorithms are expensive due to a large number of parameters of the two core tensors $\tG$ and $\tH$. The proposed algorithm 
will estimate only four vectors of length $R$  per dimension whereas the core tensors $\tG$ and $\tH$ need not to be estimated.

We will first introduce an orthogonal normalisation for the block tensor deflation, then state the correctness of the proposed deflation scheme.

\begin{lemma}[Orthogonal normalization for rank splitting]\label{lem_orthofactor2}
Given a decomposition of $\tY$ as $ \tY \approx \llbracket
\tG;\bU^{(1)}, \bU^{(2)},  \ldots, \bU^{(N)} \rrbracket + 
\llbracket \tH;  \bV^{(1)}, \bV^{(2)}, \ldots, \bV^{(N)}
\rrbracket$, where $\bU^{(n)} \in \Real^{I_n \times (K)}$ and $\bV^{(n)} \in \Real^{I_n \times (R-K)}$, $K \le R-K$, one can construct an equivalent decomposition, denoted by tildas, which has the same approximation error,  {{such that}}
\begin{itemize}
\item $	 \llbracket
\tG;\{\bU^{(n)}\} \rrbracket = \llbracket
\widetilde{\tG};\{\widetilde\bU^{(n)}\} \rrbracket$,
$	 \llbracket
\tH; \{\bV^{(n)} \} \rrbracket= \llbracket
\widetilde{\tH}; \{ {\widetilde{\bV}}^{(n)}\}\rrbracket$
\item $\widetilde\bU^{(n)}$ and ${\widetilde{\bV}}^{(n)}$ are orthogonal, i.e.,
$(\widetilde{\bU}^{(n)})^T\, \widetilde\bU^{(n)} = \bI_{K}$ and $(\widetilde{\bV}^{(n)})^T\, {\widetilde{\bV}}^{(n)} = \bI_{R-K}$.
\item and obey conditions 
$(\widetilde{\bU}^{(n)})^T \,  {\widetilde{\bV}^{(n)}}  = \left[\diag\{\bsigma_n\}, \0_{R-2K}\right]$ where $\bsigma_n = [\sigma_{n,1}, \ldots, \sigma_{n,K}] \in \Real^{K}$  and $0 \le  \sigma_{n,r} < 1$.
\end{itemize}
\end{lemma}
\begin{proof} See Appendix~\ref{sec::proof_lem_1}.
\end{proof}

\begin{theorem}[Rank splitting]{\label{theo_rank_splitting}}
A rank-$R$ tensor $\tY = \llbracket \bbeta;  \{\bB^{(n)} \}\rrbracket$ has an exact decomposition  {
{as in (\ref{equ_btd})}} $$\tY = \llbracket \tG ; \bU^{(1)}, \ldots, \bU^{(N)}\rrbracket  + \llbracket \tH ; \bV^{(1)}, \ldots, \bV^{(N)}\rrbracket$$ where {
{$\bU^{(n)} \in \Real^{I_n \times K}$ and $\bV^{(n)} \in \Real^{I_n \times (R-K)}$, $K \le R-K$}} and
\begin{itemize}
\item at least two factor matrices $\bB^{(n)} \in \Real^{I_n \times R}$ are of full column rank,
\item $\tG$ has multilinear rank-$(K,\ldots,K)$.
\end{itemize}
Then $\tG$ is a tensor of rank-$K$  and $\tH$ of rank $(R-K)$.
\end{theorem}
\begin{proof} See Appendix~\ref{sec::proof_theorem1}.
\end{proof}

\section{Alternating Subspace Update Algorithm}\label{sec::asu_blk2}

In this section, we consider order-3 tensors of size $R \times  R \times R$. Tensors of larger and unequal sizes should be compressed to this size using the Tucker decomposition\cite{Lathauwer_HOOI,Comon09,Phan_CrNc}.
We will develop an algorithm for the block tensor deflation which reduces the rank by $K = 2$. For this particular case, the core tensor $\tG$ is size of $2 \times 2 \times 2$, and the core tensor $\tH$ of size $(R-2) \times (R-2 \times (R-2)$. The factor matrices $\bU^{(n)}$ and $\bV^{(n)}$  are of size $R \times 2$ and $R \times (R-2)$, respectively.
The rank-2 block deflation has an advantage over the rank-1 tensor deflation when factor matrices have two nearly collinear components.

We denote matrices $\bar{\bV}^{(n)}  = [{\bv}^{(n)}_1, {\bv}^{(n)}_2]$ which comprise the first two columns of $\bV^{(n)}$, 
and perform reparameterization of $\bU^{(n)}$ as 
\be
\bU^{(n)} = \bW^{(n)} \diag(\bxi_n)  + {{\bar{\bV}}}^{(n)} \diag(\bsigma_n)  \, ,  \label{equ_U_WV}
\ee 
where $\bxi_n = [\xi_{n1}, \xi_{n2}]^T$, $\xi_{nr} = \sqrt{1 - \sigma_{nr}^2}$, and $\bW^{(n)} = [\bw^{(n)}_1, \bw^{(n)}_2]$ of size $R \times 2$.
$[\bW^{(n)}, \bV^{(n)}]$ are orthonormal matrices of size $R \times R$, i.e., $[\bW^{(n)}, \bV^{(n)}]^T [\bW^{(n)}, \bV^{(n)}] = \bI_R$.

Consider the following criterion to be minimized,
\be
D = \frac{1}{2} \| \tY  - \tG \times \{\bU^{(n)} \}  - \tH \times \{\bV^{(n)}\}\|_F^2  \, .\label{equ_costD}
\ee
The ALS algorithm \cite{Lath-Nion-BCM3} and the non-linear least squares  (NLS) algorithm \cite{Sorber_fusion2013} consider the same optimisation criteria.
 We will later simplify the objective function in (\ref{equ_costD}) by replacing the core tensors by their closed-form expressions and applying the above reparameterization. The objective function will finally depend only on $\bW^{(n)}$, ${{\bar{\bV}}}^{(n)}$ and $\bsigma_n$ for $n = 1, 2, 3$.

\subsection{Closed-form expressions for the core tensors}\label{sec::core_GH}

The first derivatives of the cost function $D$ in (\ref{equ_costD}) with respect to the core tensors $\tG$ and $\tH$ are given by\ 
\be
	\frac{\partial D}{\partial \tG} &=& - \tY \times \{\bU^{(n) T}\} + \tG + \bar{\tH} \times \{\diag(\sigma_n)\}  \, , \label{equ_dG}\\
	\frac{\partial D}{\partial \tH} &=& - \tY \times \{\bV^{(n) T}\} + \tG \times \left\{\left[\begin{array}{c}\diag(\bsigma_n) \\ \0_{(R-2) \times 2}\end{array}\right]\right\} + {\tH}  \,, \label{equ_dH}
\ee
where $\bar{\tH} = \tH(\text{1:2},\text{1:2},\text{1:2})$.
We obtain closed-form expressions for $\tH$ and $\tG$ as
\be
\tH &=& \tY \times \{\bV^{(n) T}\} - \tG \times \left\{\left[\begin{array}{c}\diag(\bsigma_n) \\ \0_{(R-2) \times 2}\end{array}\right]\right\}  \,  ,\label{equ_closedform_H}  \\
\tG &=& \left(\tY \times \{\bU^{(n) T}\}   - \left( \tY \times \{{\bar{\bV}}^{(n) T}\}\right)  \* \tS \right) \oslash \left(1 - \tS \* \tS \right) \, ,  \label{equ_G}
\ee
where $\tS = \bsigma_1 \circ \bsigma_2 \circ  \bsigma_3$ is a rank-1 tensor of size $2 \times 2 \times 2$, $\*$ and $\oslash$  represent the Hadamard (element-wise) product and division, respectively.

We replace $\tH$ in the cost function (\ref{equ_costD}) by its closed-form in (\ref{equ_closedform_H}), and rewrite $D$ as
\be
D &=& \frac{1}{2} \|\tY -   \tY \times \left\{\bV^{(n)} \bV^{(n) T} \right\}  - \tG \times \{\bU^{(n)}\}  +  \tG \times \{{{\bar{\bV}}}^{(n)} \diag(\bsigma_n)\} \|_F^2 \notag\\
&=& \frac{1}{2}\left( \| \tY -   \tY \times \left\{\bV^{(n)} \bV^{(n) T} \right\}  \|_F^2  + \|\tG\|_F^2 + \|\tG \times \{\diag(\bsigma_n)\}\|_F^2   \right. \notag \\
&&    - 2  \langle  \tG \times \{\bU^{(n)}\}  ,  \tG \times \{{{\bar{\bV}} \diag(\bsigma_n) }^{(n)}\}    \rangle 
\left. -  2  \langle \tY -   \tY \times \left\{\bV^{(n)} \bV^{(n) T} \right\}   ,  \tG \times \{\bU^{(n)}\} \rangle \right) \notag\\
&=& \frac{1}{2}\left( \|\tY\|_F^2 -  \|\tY \times \left\{\bV^{(n)} \bV^{(n) T} \right\}  \|_F^2 -  \langle  \tG \* (1 - \tS \* \tS) ,  \tG \rangle \right)  \, .  \label{equ_costD2}
\ee
For an index $n\in \{1,2,3\}$, define $n_1$ and $n_2$ with $n_1<n_2$ as its complement in $\{1,2,3\}$, i.e., $\{n,n_1,n_2\}=\{1,2,3\}$. 
Put
\be
	\bt^{(n)}_{r,s} &=& \tY \times_{n_1}  \bu^{(n_1)  T}_{r}    \times_{n_2}   \bu^{(n_2) T}_{s}  \;, \; \label{equ_tn}\\
	\bz^{(n)}_{r,s} &=& \tY \times_{n_1}   \bv^{(n_1) T}_{r}   \times_{n_2}     \bv^{(n_2) T}_{s}  \; , \;\;\label{equ_zn}\\
	\bd^{(n)}_{r,s} &=& \bt^{(n)}_{r,s} - \bz^{(n)}_{r,s} \;    \sigma_{n_1,r}  \,\sigma_{n_2,s}\;.   \label{equ_dn}
\ee
%
The objective function in (\ref{equ_costD2}) can be expressed as
\be
D &=& \frac{1}{2}\left( \|\tY\|_F^2 -  \|\tY \times \left\{\bV^{(n)} \bV^{(n) T} \right\}  \|_F^2   
 -  \sum_{r_1 = 1}^{2}  \sum_{r_2 = 1}^{2} \sum_{r_3 = 1}^{2} \frac{(\bu_{r_1}^{(1) T}  \bt^{(1)}_{r_2,r_3}  -  \bv^{(1) T}_{r_1}  \bz^{(1)}_{r_2,r_3} )^2}{1  -   \sigma_{1,r_1}^2   \sigma_{2,r_2}^2 \sigma_{3,r_3}^2} \right)   \notag \\
&=& \frac{1}{2}\left( \|\tY\|_F^2 -  \|\tY \times \left\{\bV^{(n)} \bV^{(n) T} \right\}  \|_F^2 
 -  \sum_{r_1, r_2, r_3 = 1}^{2} \frac{(\xi_{1,r_1} \bw_{r_1}^{(1) T}  \bt^{(1)}_{r_2,r_3}  + \sigma_{1,r_1} \bv^{(1) T}_{r_1}  \bd^{(1)}_{r_2,r_3} )^2}{1 -  \sigma_{1,r_1}^2   \sigma_{2,r_2}^2 \sigma_{3,r_3}^2} \right)  \, .\quad\label{equ_cost_D3}
\ee

\subsection{Estimation of $\bsigma_n$}\label{sec::estimate_sigma_n}

We begin with deriving update rules for $\bsigma_1 = [\sigma_{1,1}, \sigma_{1,2}]$.
As shown in the cost function in (\ref{equ_cost_D3}), the parameters $\bsigma_1$ involve only the third term. In order to estimate $\bsigma_1$, we keep other parameters fixed. Then minimization of the cost function (\ref{equ_cost_D3}) leads to maximization of the function of $\bsigma_1$
\be
\max_{\sigma_{1,1}, \sigma_{1,2}} \quad \sum_{r_1 = 1}^{2}  \sum_{r_2 = 1}^{2} \sum_{r_3 = 1}^{2} \frac{(\xi_{1,r_1} \bw_{r_1}^{(1) T}  \bt^{(1)}_{r_2,r_3}  + \sigma_{1,r_1} \bv^{(1) T}_{r_1}  \bd^{(1)}_{r_2,r_3} )^2}{1 -  \sigma_{1,r_1}^2   \sigma_{2,r_2}^2 \sigma_{3,r_3}^2} \, .
\ee
Each $\sigma_{1,r_1}$ is found as
$\sigma_{1,{r_1}} = 1/{\sqrt{1 + x_{r_1}^2}}$
where $x_{r_1}$ is solution to the problem
\be
x_{r_1} = \arg\max_{x} \, \sum_{r_2 = 1}^{2} \sum_{r_3 = 1}^{2} \frac{(\alpha_{r_2,r_3}  \, x +  \beta_{r_2,r_3} )^2}{ x^2 + 1 -  \sigma_{2,r_2}^2 \sigma_{3,r_3}^2 }  \label{equ_optimize_x_r}
\ee
$\alpha_{r_2, r_3} = \bw_{r_1}^{(1) T}  \bt^{(1)}_{r_2,r_3}$ and
$\beta_{r_2, r_3} = \bv^{(1) T}_{r_1}  \bd^{(1)}_{r_2,r_3}$. The optimal $x_{r_1}$ is a root of a polynomial of degree-8.
The other  $\sigma_{n,r}$ can be estimated similarly.
%
%

\subsection{Estimation of orthogonal components $\bW^{(n)}$ and $\bV^{(n)}$}

This section will present update rules which preserve orthogonality constrains on $\bW^{(n)}$ and $\bV^{(n)}$. Indeed we only need to update $\bW^{(n)}$ and the first two column vectors ${{\bar{\bV}}}^{(n)} = [\bv^{(n)}_1, \bv^{(n)}_2]$, whereas the last $(R-4)$ columns $[\bv^{(n)}_3, \ldots, \bv^{(n)}_{R-2}]$ are chosen as arbitrary orthogonal complement to $[\bW^{(n)}, {{\bar{\bV}}}^{(n)}]$.

Since $\bV^{(n)} \bV^{(n) T } = \bI_R - \bW^{(n)} \bW^{(n) T}$, we have
\be
\|\tY \times \left\{\bV^{(n)} \bV^{(n) T} \right\}  \|_F^2   = \tr(\bPhi_n)    - \tr(\bW^{(n) T}  \bPhi_n \bW^{(n)})
\ee
where $\bPhi_n = \bY_{(n)}  \left(\bigotimes_{k \neq n} \bV^{(n)} \bV^{(k) T} \right) \bY_{(n)}^T$ are matrices of size $R \times R$.
The cost function in (\ref{equ_cost_D3}) is rewritten as  
\be
D &=& \frac{1}{2} \left(
\|\tY\|_F^2 -  \tr(\bPhi_n) 
+ \sum_{r = 1}^{2}  \,  \bw^{(n) T}_r   \, \bQ_{n,r} \,   \bw^{(n)}_r  - \bv^{(n) T}_r   \, \bF_{n,r} \,   \bv^{(n)}_r  -  2 \bw^{(n) T}_r   \, \bK_{n,r} \,   \bv^{(n)}_r
\right) \qquad \notag \label{equ_costD4}
\ee
where
\be
	\bQ_{n,r} &=& \bPhi_n -   \xi_{n,r}^2 \sum_{k,l}  \frac{\bt^{(n)}_{k,l}  \; {\bt^{(n)}_{k,l}}^T }{1 - \sigma_{n,r}^2 \sigma_{n_1,k}^2 \sigma_{n_2,l}^2} \, ,  \\
	\bF_{n,r} &=& \sigma_{n,r}^2   \sum_{k,l}  \frac{\bd^{(n)}_{k,l}  \; {\bd^{(n)}_{k,l}}^T }{1 - \sigma_{n,r}^2 \sigma_{n_1,k}^2 \sigma_{n_2,l}^2}\, , \\
	\bK_{n,r} &=& \xi_{n,r}\sigma_{n,r}  \sum_{k,l}  \frac{\bt^{(n)}_{k,l}  \; {\bd^{(n)}_{k,l}}^T }{1 - \sigma_{n,r}^2 \sigma_{n_1,k}^2 \sigma_{n_2,l}^2}\, .
\ee
It follows that $\bW^{(n)}$ and ${{\bar{\bV}}}^{(n)}$  are solutions to the following quadratic optimisation
\be
\min  f({\bW^{(n)}, {{\bar{\bV}}}^{(n)}}) &=&  \frac{1}{2}\left(\sum_{r = 1}^{2}  \,  \bw^{(n) T}_r   \, \bQ_{n,r} \,   \bw^{(n)}_r  - \bv^{(n) T}_r   \, \bF_{n,r} \,   \bv^{(n)}_r
-2  \sum_{r = 1}^{2}  \bw^{(n) T}_r   \, \bK_{n,r} \,   \bv^{(n)}_r \right) \\
\text{subject to} && [\bW^{(n)} \, {{\bar{\bV}}}^{(n)}]^T [\bW^{(n)} \, {{\bar{\bV}}}^{(n)}] = \bI_4. \notag
\ee
Following the Crank-Nicholson-like scheme \cite{Zaiwen_2012}, we can update the orthogonal matrices $\bX_n = [\bW^{(n)}, {\bar{\bV}}^{(n)}]$ with $\bX_n^T \bX_n = \bI_4$ using the following rules
\be
	\bX_n &\leftarrow& \bX_n  - 2\tau [\bG_f ,  \bX_n] \, \left(\bI_{8}  +  {\tau}\left[\begin{array}{cc}  \bX_n^T \bG_{f} & \bI_4 \\ - \bG_f^T \bG_f  & - \bG_f^T \bX_n \end{array}\right] \right)^{-1}  
	\left[\begin{array}{c}  \bI_4 \\ - \bG_f^T \bX_n \end{array}\right]
, \quad\label{equ_update_WV}
\ee
where $\bG_f = [\bg_{f,\bw^{(n)}_1}, \bg_{f,\bw^{(n)}_2}, \bg_{f,\bv^{(n)}_1}, \bg_{f,\bw^{(n)}_2}]$ of size $R \times 4$ are the first order derivatives of the function $ f({\bW^{(n)}, {{\bar{\bV}}}^{(n)}})$ with respect to $[\bW^{(n)}, \, {{\bar{\bV}}}^{(n)}]$
\be
	\bg_{f,\bw^{(n)}_r} &=& \frac{\partial f}{\partial \bw^{(n)}_r} = \bQ_{n,r} \, \bw^{(n)}_r  -  \bK_{n,r} \, \bv^{(n)}_r   \, , \label{equ_grad_f_wn} \\
	\bg_{f,\bv^{(n)}_r}  &=& \frac{\partial f}{\partial \bv^{(n)}_r} =  - \bF_{n,r} \, \bv^{(n)}_r  -  \bK_{n,r}^T \, \bw^{(n)}_r \,  \label{equ_grad_f_vn},
\ee
and $\bGamma_n = \bX_n^T \bG_f$ and $\tau > 0$ is a step size chosen using the Barzilai-Borwein method\cite{citeulike:2678977}.
Each iteration to update $\bX_n = [\bW^{(n)}, {{\bar{\bV}}}^{(n)}]$ inverts a matrice of size $4 \times 4$.

We finally derive update rules for all parameters. The proposed Alternating Subspace Update (ASU) algorithm is summarized in Algorithm~\ref{alg_asu}.
The algorithm alternating updates $\bsigma_n$ and $[\bW^{(n)}, {{\bar{\bV}}}^{(n)}]$ for $n = 1, 2, 3$. The entire factor matrices $\bV^{(n)}$ and core tensors $\tG$, $\tH$ are computed only once.

  {
\setlength{\algomargin}{1em}
\begin{algorithm}[t!]
\SetFillComment
\SetSideCommentRight
\CommentSty{\footnotesize}
\caption{{\tt{Alternating Subspace Update (ASU)}}\label{alg_asu}}
\DontPrintSemicolon \SetFillComment \SetSideCommentRight
\KwIn{Data tensor $\tY$:  $(R \times R  \times R)$ of rank $R$}
\KwOut{A rank-(2,2,2) tensor $\llbracket \tG; \{\bU^{(n)}\} \rrbracket$ and rank-$(R-2,R-2,R-2)$ tensor $\llbracket \tH; \{\bV^{(n)}\} \rrbracket$} \SetKwFunction{BCD}{BCD\_(2,L)}
\SetKwFunction{evd}{eigs}
\Begin{
\nl Initialise components $\bU^{(n)}$ and $\bV^{(n)}$\;
\nl Orthogonal normalization to $\bU^{(n)}$ and $\bV^{(n)}$ and compute $\bsigma_n = [\sigma_{n,1}, \sigma_{n,2}]^T$ and $\bW^{(n)}$\;
\Repeat{a stopping criterion is met}{
\For{$n = 1, 2, 3$}{
 \For{$r = 1, 2$}{
\nl	 	 Update $\sigma_{n,r} = \frac{1}{\sqrt{1 + x^2}}$ where $x$ is solved as in (\ref{equ_optimize_x_r})
	}
\nl Compute $\bG_{f}$ as in (\ref{equ_grad_f_wn}) and (\ref{equ_grad_f_vn}), $\bGamma_n = \bX_n^T \bG_{f}$ where $\bX_n = [\bW^{(n)}, {{\bar{\bV}}}^{(n)}]$\;
\nl Update $\bX_n = [\bW^{(n)}, {{\bar{\bV}}}^{(n)}]$ as in (\ref{equ_update_WV})\;
    %
\nl $\bU^{(n)}  \leftarrow  \bW^{(n)} \, \diag(\bxi_n) + {{\bar{\bV}}}^{(n)} \,
\diag(\bsigma_n) $ \; } }
	  \For{$n = 1, \ldots, N$}{
\nl Select $\bV^{(n)}_{3:R-2}$ as an orthogonal
complement of $[\bW^{(n)}, {{\bar{\bV}}}^{(n)}]$}
 \nl Compute output $\tG$ and $\tH$ as in (\ref{equ_G})  and (\ref{equ_closedform_H})\;
}
\end{algorithm}
\def\baselinestretch{1}}


The most expensive step in the ASU algorithm is computation of the matrices $\bPhi_n = \bY_{(n)}  \left(\bigotimes_{k \neq n} \bV^{(n)} \bV^{(k) T} \right) \bY_{(n)}^T$. A naive computation method might cost $\calO(R^4)$.
We present a more efficient computation which requires a cost of order $\calO(R^3)$
\be
	\bPhi_n &=& \bY_{(n)}  \left((\bI  - \bW^{(n_2)} \bW^{(n_2) T} ) \otimes   (\bI - \bW^{(n_1)} \bW^{(n_1) T}) \right) \bY_{(n)}^T \notag\\
	&=& \bY_{(n)} \,  \bY_{(n)}^T   -   \bY_{(n)}  (\bW^{(n_2)} \bW^{(n_2) T} \otimes \bI)  \bY_{(n)}^T   
	\notag\\  && 
	-  \bY_{(n)}  \, ( \bI  \otimes  \bW^{(n_1)} \bW^{(n_1) T})      \bY_{(n)}^T 
	+  \bY_{(n)} (\bW^{(n_2)} \bW^{(n_2) T} \otimes \bW^{(n_1)} \bW^{(n_1) T})  \bY_{(n)}^T \notag\\
	&= & \bY_{(n)} \,  \bY_{(n)}^T - \langle \tY \times_{n_1}  \bW^{(n_1)} , \tY \times_{n_1}  \bW^{(n_1)} \rangle_{n_1,n_2}   \notag \\ && - \langle \tY \times_{n_2}  \bW^{(n_2)} , \tY \times_{n_2}  \bW^{(n_2)} \rangle_{n_1,n_2}  
	- \langle \tY \times_{n_1}  \bW^{(n_1)} \times_{n_2}  \bW^{(n_2)} , \tY \times_{n_1}  \bW^{(n_1)} \times_{n_2}  \bW^{(n_2)} \rangle_{n_1,n_2} \, ,  \notag \qquad \label{equ_Phin}
\ee
where $\{n_1 < n_2\} = \{1,2,3\} \setminus \{n\}$.

The first term $ \bY_{(n)} \,  \bY_{(n)}^T$ is computed only once.
The mode-$n_k$ tensor productions $\tY \times_{n_k}  \bW^{(n_k)}$ yields a tensor comprising two slices of size $R \times R$  with a computation cost of $\calO(R^3)$. 

\section{Simulations}\label{sec::simulations}

\begin{example}[Decomposition of small tensors admitting the CP model.]\label{ex_1}
\end{example}
In this first example, we illustrate the block deflation of tensor of size $R\times R \times R$ and of rank $R$ where $R = 10, 20, 30$.
The weight coefficients $\lambda_r$ were set to 1, whereas collinearity degrees between components $\ba^{(n)}_r$ and $\ba^{(n)}_s$ for all $r \neq s$ were identical to a specific value $c$, which was varied 
in the range [0, 0.9], $\ba^{(n) T}_r \ba^{(n)}_s = c$ and $\ba^{(n) T}_r \ba^{(n)}_r  = 1$ for all $n$ (see Appendix~F in \cite{Phan_tensordeflation_init}). We use the subroutine $``\tt{gen\_matrix}''$  in the TENSORBOX\cite{Phan_Tensorbox} to generate factor matrices with specific correlation coefficients.

We compare the ASU algorithm with the ALS algorithm \cite{Lath-Nion-BCM3} for the multilinear rank-$(L_r,M_r,N_r)$ block tensor decomposition with two blocks. 
For this problem, one can use the non-linear least squares (NLS) algorithm \cite{Sorber_fusion2013}. However, as similar to the ALS algorithm \cite{Lath-Nion-BCM3}, the NLS algorithm needs to estimate two core tensors and full factor matrices. Hence this algorithm is much more expensive than the ASU algorithm. Simulations were run on a Macbook-air laptop having 4 GB memory and a 1.8 GHz core i7. Due to space and time consuming, the ALS \cite{Lath-Nion-BCM3} was only ran in simulations for $R = 10$.

The algorithms were initialised by the same values generated using the Direct Trilinear Decomposition (DTLD) \cite{Sanchez1990}.
The algorithms ran until differences between consecutive approximation errors were small enough, $|\varepsilon_{k} -  \varepsilon_{k+1}| \le 10^{-6} \,  \varepsilon_{k}$ where $\varepsilon =  \|\tY - \hat{\tY}\|_F^2$, or when the number of iterations exceeded 1000.
Rank-1 tensors were then obtained from decomposition of blocks of rank-2.
Performances were assessed through the squared angular errors SAE in estimation of components $\ba^{(n)}$ SAE $= \arccos\left(\frac{\ba^T \, \hat{\ba}}{\|\ba\|_2\|\hat{\ba}\|_2}\right)^2$.
There were 100 independent runs for each rank $R = 10, 20$ and 30. The Gaussian noise was added into the tensor with signal-noise-ratio SNR = 30 dB.

\begin{figure}[t]
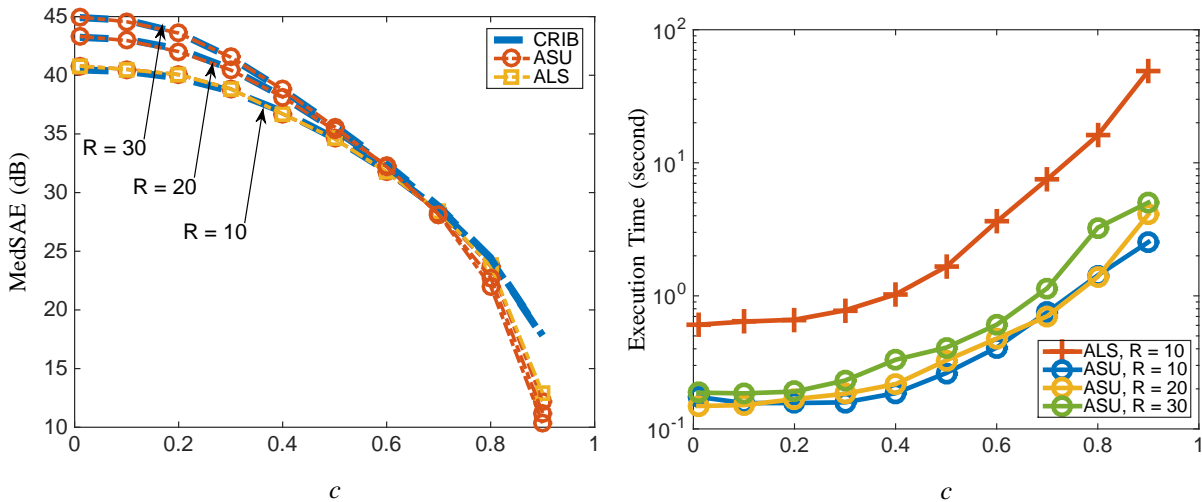

\centering
{
\psfrag{MSAE (dB)}[b][b]{\scalebox{1}{\color[rgb]{0,0,0}\setlength{\tabcolsep}{0pt}\begin{tabular}{c}\footnotesize MedSAE (dB)\end{tabular}}}%
\psfrag{c}[t][t]{\scalebox{1}{\color[rgb]{0,0,0}\setlength{\tabcolsep}{0pt}\begin{tabular}{c}\small $c$\end{tabular}}}%
\includegraphics[width=.47\linewidth, trim = 0.0cm -0.2cm 0cm 0cm,clip=true]{blkdefla_2__sae}\label{fig_ex1_sae}}
{
\psfrag{Execution Time (second)}[c][c]{\scalebox{1}{\color[rgb]{0,0,0}\setlength{\tabcolsep}{0pt}\begin{tabular}{c}\footnotesize Execution Time (second)\end{tabular}}}%
\psfrag{c}[t][t]{\scalebox{1}{\color[rgb]{0,0,0}\setlength{\tabcolsep}{0pt}\begin{tabular}{c}\small $c$\end{tabular}}}%
\includegraphics[width=.47\linewidth, trim = 0.0cm -0.2cm 0cm 0cm,clip=true]{blkdefla_2_exectime}\label{fig_ex1_rtime}}
\caption{Comparison of median SAEs and execution times of the ASU and ALS algorithms \cite{Lath-Nion-BCM3} in decomposition of tensors of size $R \times R \times R$ and rank $R$ where $R$ = 10, 20 and 30 for Example~\ref{ex_1}.}
\label{fig_sae_R10_30}
\end{figure}

%
Fig.~\ref{fig_sae_R10_30} shows median SAE (MedSAE) in dB ($-10\log_{10} SAE$) obtained by ASU and ALS\cite{Lath-Nion-BCM3} compared with the Cram{\'e}r-Rao Induced bound (CRIB) \cite{Petr_CRIB} on the squared angular error.
Algorithms succeeded in most cases, but failed only when $c = 0.9$. For such difficult scenarios, CRIB on SAE was about 17.8 dB, indicating median angular error of 7.4 degrees between the original and estimated components.
We note that in practice, it is hard to estimate a component with CRIB less than 20 dB, i.e.,  angular error of 5.7 degrees \cite{Phan_fLM}.

In Fig.~\ref{fig_ex1_rtime}, we compare execution times (in second) of algorithms for different ranks. Since the decomposition became more difficult when $c$ was close to 1, running times of algorithms increased as shown in  Fig.~\ref{fig_ex1_rtime}. The ASU algorithm was on average 8 times faster than ALS \cite{Lath-Nion-BCM3} when $R = 10$. 

The results confirmed high speed and accuracy of the proposed ASU algorithm.

\begin{example}[Decomposition of large-scale tensors with high rank]\label{ex_2}
\end{example}
This example illustrates an advantage of ASU over existing algorithms for the ordinary CPD in decomposition of large-scale tensors with relatively high rank $R$ = 300 and 500.
We generated rank-$R$ synthetic tensors of size $R \times R \times R$ as in the previous example. Components $\ba^{(n)}_r$ and $\ba^{(n)}_s$ for $r \neq s$ have identical collinearity degrees, i.e., $\ba^{(n) T}_r \ba^{(n)}_s = c$ where $c = 0.1, 0.2, \ldots, 0.6$. The Gaussian noise was at SNR = 30 dB.
Simulations were run on a computer consisted of Intel Xeon 2 processors clocked at 3.33 GHz, 64GB of main memory.
 Extraction of all components is expensive in both computation time and space. The main reason is that CP gradient computation is with a cost of $\calO(R^4)$ \cite{Phan_FastALS}.
For such big tensors, sequential extraction of rank-1 tensors using the ASU algorithm is more efficient. The ASU algorithm is particularly suited to tracking a few components without estimation of the full CP model as other algorithms.
In this example, ASU could extract components after, on average, only 3.8 seconds for $R$ = 300, and 20 seconds when $R = 500$. Decomposition of the same tensors using the FastALS algorithm for CPD \cite{Phan_FastALS} on average needed 538 and 3675 seconds, respectively.
Comparison of execution times of ASU and FastALS\cite{Phan_FastALS} is given in Table~\ref{table_exectime}.


\begin{table}[t]
\caption{Comparison of execution times of the ASU algorithm to extract two components from  high rank-$R$ tensors, and those of the CP-FastALS algorithm for Example~\ref{ex_2}.
}\label{table_exectime}
\centering
\begin{tabular}{lcccccccccc}
			&   \multicolumn{6}{c}{Execution time (second)}  \\\cline{2-7}
			&  $c$ =  0.1  & 0.2 & 0.3 & 0.4  & 0.5  & 0.6  \\\hline
$R = 300$ \\\cline{1-1}
ASU			&  3.81   &  3.66    & 3.76    & 3.82  & 3.89    & 3.77 \\
CP-FastALS	& 530.6  & 543.5  & 537.6  & 537.6  & 541.9  &539.2  \\
\hline
$R = 500$ \\\cline{1-1}
ASU			&  38.4    & 16.7     & 16.5   & 16.9  & 16.8   & 17.1   	 \\
CP-FastALS	& 3658    & 3672   &  3679    & 3693 & 3678   & 3669     \\ \hline
%
\end{tabular}
\end{table}

\begin{example}[Comparison of rank-1 and block tensor deflations]\label{ex_3}
\end{example}

\begin{figure}[t]
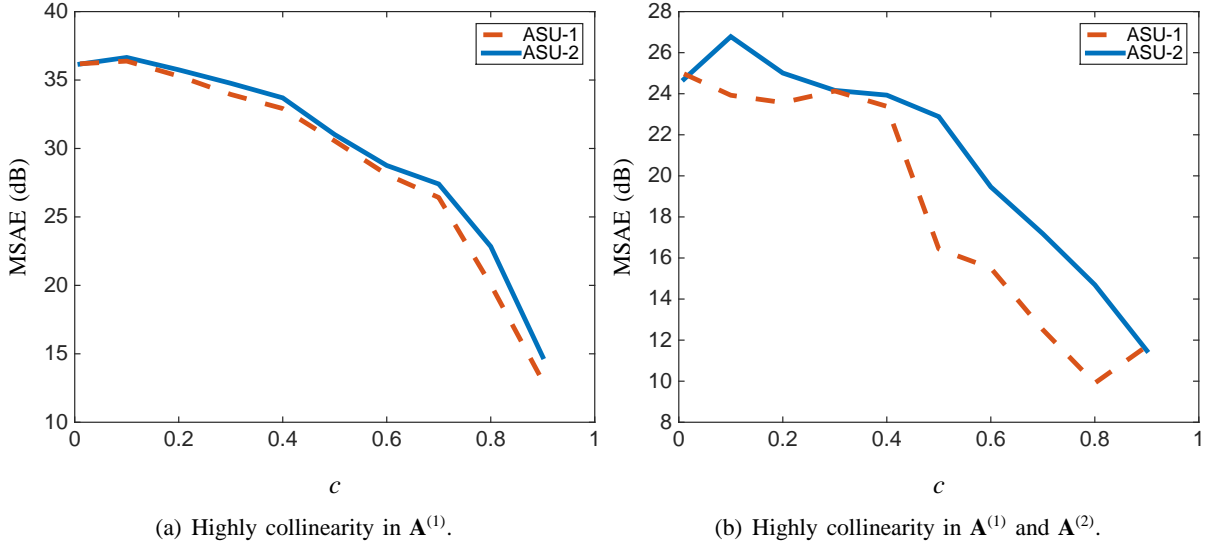

\centering
\subfigure[Highly collinearity in $\bA^{(1)}$.]
{
\psfrag{MSAE (dB)}[b][b]{\scalebox{1}{\color[rgb]{0,0,0}\setlength{\tabcolsep}{0pt}\begin{tabular}{c}\footnotesize MSAE (dB)\end{tabular}}}%
\psfrag{c}[t][t]{\scalebox{1}{\color[rgb]{0,0,0}\setlength{\tabcolsep}{0pt}\begin{tabular}{c}\small $c$\end{tabular}}}%
\includegraphics[width=.47\linewidth, trim = 0.0cm -0.2cm 0cm 0cm,clip=true]{blkdefla_2_R10_collinear_sae}\label{fig_ex3_sae}}
\subfigure[Highly collinearity in $\bA^{(1)}$ and $\bA^{(2)}$.]
{
\psfrag{MSAE (dB)}[b][b]{\scalebox{1}{\color[rgb]{0,0,0}\setlength{\tabcolsep}{0pt}\begin{tabular}{c}\footnotesize MSAE (dB)\end{tabular}}}%
\psfrag{c}[t][t]{\scalebox{1}{\color[rgb]{0,0,0}\setlength{\tabcolsep}{0pt}\begin{tabular}{c}\small $c$\end{tabular}}}%
\includegraphics[width=.47\linewidth, trim = 0.0cm -0.2cm 0cm 0cm,clip=true]{blkdefla_2_R10_2collinear_sae}\label{fig_ex3_sae_2col}}
\caption{Comparison of mean SAEs (MSAE) achieved by the ASU algorithms for rank-1 tensor deflation\cite{Phan_tensordeflation_alg} and block tensor deflation for Example~\ref{ex_3}.}
\label{fig_sae_collinear}
\end{figure}

This example presents a case when the block tensor deflation is more appropriate than the rank-1 tensor deflation. We considered tensors whose factor matrix $\bA^{(1)}$ comprised two highly collinear components.
More specifically, we first generated rank-$R$ synthetic tensors of size $R \times R \times R$ where $R = 10$ as tensors in Example~\ref{ex_1}, i.e., $\ba^{(n) T}_r  \ba^{(n)}_r = 1$ and $\ba^{(n) T}_r \ba^{(n)}_s  = c$ for all $r \neq s$ and $0<c < 1$.
The component $\ba^{(1)}_2$ was then adjusted so that its collinearity degree with $\ba^{(1)}_1$ was of $\rho = 0.98$
\be
	\ba^{(1)}_2 :=   (\rho - c \alpha) \, \ba^{(1)}_1 + \alpha  \ba^{(1)}_2	
\ee
where $\alpha =  \sqrt{(1 - \rho^2)/(1-c^2)}$. Collinearity degrees between $\ba^{(1)}_2$ and the other components $\ba^{(1)}_r$ for $r > 2$ were then given by
\be
\ba^{(1) T}_2  \, \ba^{(1)}_r = c (\rho +  \alpha (1-c)) \, .
\ee
Since $\ba^{(1)}_1$ or $\ba^{(1)}_2$ were highly collinear, extraction of only one rank-1 tensor associated with $\ba^{(1)}_1$ or $\ba^{(1)}_2$ is difficult  as analysed in Part 2 \cite{Phan_tensordeflation_init}.
We will show that there are loss of accuracy in extraction of the rank-1 tensor $\ba^{(1)}_1 \circ \ba^{(2)}_1 \circ \ba^{(3)}_1$, compared with block tensor deflation which extracts two rank-1 tensors comprising components $\ba^{(1)}_1$ or $\ba^{(1)}_2$.
For this comparison, we initialised the ASU algorithm (ASU-1) \cite{Phan_tensordeflation_alg} for the rank-1 tensor deflation and the ASU algorithm proposed in this paper (ASU-2) by the true components.
The mean SAEs (dB) of estimated components achieved by the two algorithms shown in Fig.~\ref{fig_sae_collinear} indicate that the loss varied from 0.37 dB to 2.5 dB when $c$ increased from 0.1 to 0.9.

In another simulation with similar settings, we compared ASU-1 and ASU-2 when the factor matrices $\bA^{(1)}$ and $\bA^{(2)}$ comprised two highly collinear components $\ba^{(1) T}_1 \ba^{(1)}_2  =  \ba^{(2) T}_1 \ba^{(2)}_2 = 0.95$. It is necessary to remind conditions for the rank-1 tensor deflation, i.e, conditions for ASU-1.
According to Lemma 2 in Part 1\cite{Phan_tensordeflation_alg}, a rank-1 tensor can only be uniquely extracted if at least two components do not lie 
within the column spaces of the other components. 
Since the two components $\ba^{(1)}_1$ and $\ba^{(2)}_1$ were highly collinear with $\ba^{(1)}_2$ and  $\ba^{(2)}_2$, respectively, the rank-1 tensors $\ba^{(1)}_1 \circ \ba^{(2)}_1 \circ \ba^{(3)}_1$ and $\ba^{(1)}_2 \circ \ba^{(2)}_2 \circ \ba^{(3)}_2$ can be considered to violate the condition. 
Extraction of one of the two rank-1 tensors is not stable. Instead, they should be extracted together. It is shown in Fig.~\ref{fig_ex3_sae_2col} that the loss of accuracy of ASU-1 was higher for this difficult decomposition. 

\section{Conclusions}\label{sec::conclusion}

We have introduced a rank-splitting scheme for CPD, and developed an ASU algorithm for rank-2 block deflation. The algorithm needs to estimate only 4 vectors and two scalars per dimension, and has a computational cost of $\calO(R^3)$ for a tensor of size $R \times R \times R$.
The algorithm can be extended to higher order tensors, and decomposition with additional constraints. 
Algorithms for the block tensor deflation are implemented in the Matlab package TENSORBOX which is available online at: {\url{http://www.bsp.brain.riken.jp/~phan/tensorbox.php}}.

\appendices
\section{Proof of Lemma~\ref{lem_orthofactor2}}\label{sec::proof_lem_1}

{
{
\begin{proof}  Let $\bQ_n$ and $\bF_n$ be column space of $\bU^{(n)}$, and
$\bV^{(n)}$, respectively, which can be obtained from QR decompositions
\be
	\bU^{(n)} = \bQ_{n} \, \bR_{n}  \, ,  \qquad
	\bV^{(n)} = \bF_{n} \, \bK_{n}  \, . \notag
\ee
Consider singular value decomposition (SVD) of  ${\bQ_n}^T \,
\bF_{n}  =  \bGamma_n  \, \bSigma_n \, \bPsi_n^T$ where $\bGamma_n \; \in   \Real^{K \times K}$, $\bPsi_n \; \in   \Real^{K \times (R-K)}$ and $\bSigma_n = \left[\diag\{\bsigma_n\} , \0_{R-2K} \right]$, $\bsigma_n \in \Real_+^{K}$.
Then, the new decomposition is equivalently defined through
\begin{eqnarray}
    \widetilde\bU^{(n)} &=& \bQ_n \, \bGamma_n,\qquad n=1,\ldots,N, \\
    \widetilde{\tG} &=&   \tG \times_1 \, (\bGamma_1^T\, \bQ_1^T \bU^{(1)}) \,  \cdots \times_N \, (\bGamma_N^T \,  \bQ_N^T \, \bU^{(N)})\, ~,
\end{eqnarray}
and
\begin{eqnarray}
    {\widetilde{\bV}}^{(n)} &=& \bF_n \, \bPsi_n,\qquad n=1,\ldots,N, \\
    \widetilde{\tH} &=&   \tH \times_1 \, (\bPsi_1^T\, \bF_1^T \bV^{(1)}) \,  \cdots \times_N \, (\bPsi_N^T \,  \bF_N^T \, \bV^{(N)})~.
\end{eqnarray}
It can be verified that $\widetilde\bU^{(n)}$ and ${\widetilde{\bV}}^{(n)}$ are orthogonal and
\be
(\widetilde\bU^{(n)})^T \, {\widetilde{\bV}}^{(n)} = \bGamma_n^T \,
\bQ_n^T \, \bF_{n}\, \bPsi_{n} = \bSigma_n. \ee
This completes the proof.
\end{proof}}

\section{Proof of Theorem~\ref{theo_rank_splitting}}\label{sec::proof_theorem1}
\begin{proof}
For simplicity, we assume that $\bB^{(1)}$ and $\bB^{(N)}$ are of full column rank. Since
\be
	\bY_{(n)} = \bB_{n}  \diag\{\bbeta\} \left( \bigodot_{k \neq n} \, \bB^{(k)}  \right)
	 = \left[\bU^{(n)} \, \bV^{(n)}  \right] \, \left[\begin{array}{@{}c@{}}
	 \bG_{(n)} \, \left(\bigodot_{k \neq n} \, \bU^{(k)} \right)^T \\
	 \bH_{(n)} \, \left(\bigodot_{k \neq n} \, \bV^{(k)} \right)^T  \\	
	 \end{array}\right]\, , \notag
\ee
$\bU^{(1)}$, $\bV^{(1)}$ and $\bU^{(N)}$, $\bV^{(N)}$ are also full column rank matrices.

Thanks to Lemma~\ref{lem_orthofactor2}, we can assume, without any loss in generality, that the factor matrices $\bU^{(n)}$ and $\bV^{(n)}$ for $n = 1$ and $n = N$, obey the normalization condition, i.e., ${\bU^{(n)}}^T \,  \bU^{(n)}  = \bI_{K}$, ${\bV^{(n)}}^T \,  \bV^{(n)}  = \bI_{R-K}$ and ${\bU^{(n)}}^T \,  \bV^{(n)}   = \left[ \diag(\bsigma_n), \0_{K \times( R-2K)} \right]$ where $\bsigma_n  = [\sigma_{n,1}, \sigma_{n,2}, \ldots, \sigma_{n,K}]^T \in \Real^{K}$, and $0 \le  \sigma_{n,k} < 1$.

Let ${\bZ}_{N} = \left[{\bz}^{(N)}_1 , \ldots, {\bz}^{(N)}_K\right]$ be an $I_N \times K$ matrix whose columns are defined as
\be
\displaystyle {\bz}^{(N)}_k = \frac{\bu^{(N)}_k - \sigma_{N,k} \, \bv^{(N)}_{k}}{1-\sigma_{N,k}^2}, \quad  k = 1, \ldots, K.
\ee
We have ${\bZ}_{N}^T \,  \bV^{(N)}  = \0$ and ${\bZ}_{N}^T \, {\bU}^{(N)} = \bI_{K}$. Put $\bW = {\bZ}_{N}^T \, \bB^{(N)}$, the tensor-matrix product $\tY \, {\times}_N \, {\bZ}_{N}^T$ is given by
\be	
	\tY \, {\times}_N \, {\bZ}_{N}^T
	&=& \llbracket  \bbeta_{\calR} \, ;  \bB_{\calR}^{(1)}, \ldots, \bB_{\calR}^{(N-1)},  \bW_{\calR} \rrbracket   , \label{equ_Ytu1_a}
\ee
where $\calR$ denotes set of indices of non-zero columns $\bw_{k} \neq 0$ for $k \in \calR$, $\bB^{(n)}_{\calR} = \bB^{(n)}(:,{\calR})$ are sub matrices taken from $\bB^{(n)}$ and $\bbeta_{\calR} = \bbeta(\calR)$. 

From the block term decomposition of $\tY$, we also have
\be	
	\tY \, {\times}_N \, \bZ_{N}^T  =  \llbracket \tG ; \bU^{(1)}, \ldots, \bU^{(N-1)}  ,  \bI_K \rrbracket \, , \label{equ_Ytu1_c}
\ee
which leads to 
\be
\tG = \llbracket  \bbeta_{\calR} \, ;  \bU^{(1) T} \, \bB_{\calR}^{(1)}, \ldots, \bU^{(N-1) T} \, \bB_{\calR}^{(N-1)},  \bW_{\calR} \rrbracket  \, . \label{equ_Ytu1_e}
\ee
Hence, the expression in (\ref{equ_Ytu1_c}) is equivalently rewritten as 
\be
\tY \, {\times}_N \, \bZ_{N}^T  =  \llbracket \bbeta_{\calR} \, ;  \bU^{(1)} \bU^{(1) T} \, \bB_{\calR}^{(1)}, \ldots, \bU^{(N-1)} \bU^{(N-1) T} \, \bB_{\calR}^{(N-1)},  \bW_{\calR} \rrbracket \,\label{equ_Ytu1_d}
\ee
Since $\bB^{(1)}_{\calR}$ is a full-column rank matrix, the CPDs in  (\ref{equ_Ytu1_a}) and (\ref{equ_Ytu1_d}) are unique and therefore identical. It follows that 
\be
	(\bI_K -  \bU^{(n)} \bU^{(n) T}) \, \bB_{\calR}^{(n)}  = 0 , \quad n = 1, \ldots, N-1.
\ee
That is $\bB_{\calR}^{(n)}$ are spanned by $\bU^{(n)}$ for $n = 1, \ldots, N-1$, respectively.
In addition, since $\tG$ has multilinear rank-$(K,\ldots,K)$, from (\ref{equ_Ytu1_e}), $\bB_{\calR}^{(1)}$ must be of size $I_1 \times K$, and can be expressed as
\be
	\bB_{\calR}^{(1)} = \bU^{(1)} \, \bQ_1
\ee
where $\bQ_1$ is a full-column rank matrix of size $K \times K$.
Implying that $\tG$ is a rank-$K$ tensor, and uniquely identified
\be
	\tG =  \llbracket  \bbeta_{\calR} \, ;  \bQ_1, \bU^{(2) T} \, \bB_{\calR}^{(2)}, \ldots, \bU^{(N-1) T} \, \bB_{\calR}^{(N-1)},  \bW_{\calR} \rrbracket  \, .  \label{equ_G_N}
\ee

Similarly we can prove that
\be
\tY \, {\times}_1 \, \bZ_{1}^T  =  \llbracket \bbeta_{\calX} \, ;  \bU^{(1)} \bZ_1^{T} \, \bB_{\calX}^{(1)}, \bU^{(2)} \bU^{(2) T} \, \bB_{\calX}^{(2)}, \ldots,  {\bU}^{(N)}  {{\bU}^{(N) T}} \bB^{(N)}_{\calX} \rrbracket \,\notag
\ee
and 
%
\be
	\tG &=& \llbracket  \bbeta_{\calX} \, ; {{\bZ}_1^T} \bB^{(1)}_{\calX} , \ldots, \bU^{(N-1) T} \, \bB_{\calX}^{(N-1)},  {{\bU}^{(N) T}} \bB^{(N)}_{\calX} \rrbracket  \, , \quad \label{equ_G_1}
\ee
where $\calX$ is an index set of $K$ non-zero columns $\bZ_1^T \bB^{(1)}$. 

Since the first and the last factor matrices in the CP decompositions of G in (\ref{equ_G_N}) and in (\ref{equ_G_1}) are of full column rank, the decompositions are unique. Therefore, the two sets $\calR$ and $\calX$ are identical, and the tensor $\llbracket \tG ; \bU^{(1)}, \ldots, \bU^{(N)} \rrbracket$ is a rank-$K$ tensor taken from $K$ rank-1 tensors of the tensor $\tY$, 
\be
	\llbracket \tG ; \bU^{(1)}, \ldots, \bU^{(N-1)}, \bU^{(N)} \rrbracket = \llbracket  \bbeta_{\calR} \, ;  \bB_{\calR}^{(1)}, \ldots, \bB_{\calR}^{(N-1)}, \bB_{\calR}^{(N)}\rrbracket \,.
\ee
Finally, it is obvious that eliminating the rank-$K$ tensor $\llbracket \tG ; \bU^{(1)}, \ldots, \bU^{(N)} \rrbracket$ from $\tY$ remains a rank-$(R-K)$ tensor, i.e. $\tH$ is a rank-$(R-K)$ tensor.

\end{proof}

\bibliographystyle{IEEEbib}

\end{document}